\documentclass[a4paper,12pt]{article}
\usepackage{latexsym}
\usepackage{amssymb}
\usepackage{amsfonts}
\usepackage{amsmath,amsthm}
\usepackage{mathrsfs}
\usepackage{bm}
\usepackage{enumerate}
\usepackage{stmaryrd}
\usepackage{xcolor}
\usepackage{commath}
\usepackage{array}
\usepackage{cases}
\usepackage{comment}

\newenvironment{@abssec}[1]{%
    \if@twocolumn

      \section*{#1}%
    \else

      \vspace{.05in}\footnotesize
      \parindent .2in
 {\upshape\bfseries #1. }\ignorespaces
    \fi}

    {\if@twocolumn\else\par\vspace{.1in}\fi}

\newenvironment{keywords}{\begin{@abssec}{\keywordsname}}{\end{@abssec}}

\newenvironment{AMS}{\begin{@abssec}{\AMSname}}{\end{@abssec}}

\newcommand\keywordsname{Key words}
\newcommand\AMSname{AMS subject classifications}
\newcommand\AMname{AMS subject classification}
\newcommand\restr[2]{{
\left.\kern-\nulldelimiterspace 
#1 
\vphantom{|} 
\right|_{#2} 
}}

\newtheorem{thm}{Theorem}[section]
\newtheorem{lem}[thm]{Lemma}

\newtheorem{prop}[thm]{Proposition}
\newtheorem{rem}[thm]{Remark}
\newtheorem{dfn}[thm]{Definition}

\newtheorem{e.g.}{Example}

\def\XXint#1#2#3{{\setbox0=\hbox{$#1{#2#3}{\int}$}
\vcenter{\hbox{$#2#3$}}\kern-.5\wd0}}

\newcommand{\R}{\mathbb R}
\newcommand{\Nt}{\mathbb N}

\newcommand{\kp}{\kappa}

\newcommand{\Gm}{\Gamma}

\newcommand{\ep}{\varepsilon}

\newcommand{\s}{\sigma}

\def\<{\langle }

\newcommand{\dl}{\delta}
\newcommand{\bt}{\beta}

\newcommand{\eq}[1]{\begin{equation}#1\end{equation}}

\newcommand{\eqn}[1]{\begin{equation*}#1\end{equation*}}

\allowdisplaybreaks[1]

\title{\bf Convergence rate of Tsallis entropic regularized optimal transport}

\author{Takeshi Suguro and Toshiaki Yachimura}
\date{}

\begin{document}

\maketitle

\begin{abstract}
In this paper, we study the Tsallis entropic regularized optimal transport in the continuous setting and establish fundamental results such as the $\Gamma$-convergence of the Tsallis regularized optimal transport to the Monge--Kantorovich problem as the regularization parameter tends to zero. In addition, using the quantization and shadow arguments developed by Eckstein--Nutz, we derive the convergence rate of the Tsallis entropic regularization and provide explicit constants. Furthermore, we compare these results with the well-known case of the Kullback--Leibler (KL) divergence regularization and show that the KL regularization achieves the fastest convergence rate in the Tsallis framework.
\end{abstract}



\pagestyle{plain}
\thispagestyle{plain}

\section{Introduction}\label{intro}
Optimal transport theory is widely used to measure distances between probability measures and to compare them. Optimal transport has been developing rapidly in recent years. In particular, the regularization of optimal transport has made it possible to compute even for high-dimensional data. 

Let $(X_1,\mu_1)$ and $(X_2,\mu_2)$ be Polish probability spaces and $c: X_1 \times X_2 \to [0, +\infty)$ be a continuous cost function. The following regularized optimal transport is the most commonly used: 
\eq{\label{KL reg OT}
  \mathrm{OT}_{\ep} = \inf_{\pi \in \Pi(\mu_1, \mu_2)} \left\{\int_{X_1 \times X_2} c\, d\pi + \ep \, \mathrm{KL}(\pi, \mu_1 \otimes \mu_2)\right\},  
}
where $\Pi(\mu_1, \mu_2)$ is the set of couplings of the given marginals $\mu_1$ and $\mu_2$. Also, $\ep > 0$ is the regularization parameter and $\mathrm{KL}$ denotes the Kullback--Leibler divergence \cite{KL51} of the transport plan $\pi$ with respect to the product measure $\mu_1 \otimes \mu_2$, which is defined for general probability measures as
\eq{\label{KL divergence}
  \mathrm{KL}(\mu, \nu) = \left\{
  \begin{aligned}
  &\int \frac{d\mu}{d\nu} \log{\frac{d\mu}{d\nu}}\, d\nu
    && \mathrm{if} \ \mu \ll \nu, \\
  &\infty 
    && \mathrm{otherwise}.
  \end{aligned}
  \right.
}
In the case where $\ep = 0$, we note that the problem \eqref{KL reg OT} consistent with the following Monge--Kantorovich problem (we refer to the standard books \cite{Ambrosio2021, Figalli2021, Santambrogio, villani2009optimal, villani2021topics} and references therein): 
\eq{\label{MK pb}
  \mathrm{OT} = \inf_{\pi \in \Pi(\mu_1, \mu_2)} \int_{X_1 \times X_2} c\, d\pi. 
}

The seminal work first introduced entropic regularization for the optimal transport problem by one of the most famous physicists E. Schr\"{o}dinger in the viewpoint of statistical physics \cite{schrodinger1931}. 
Schr\"{o}dinger studied the following minimization problem: 
\begin{equation}\label{Schrodinger bridge pb}
\inf_{\pi \in \Pi(\mu_1, \mu_2)} \ep \mathrm{KL}(\pi, \mathcal{K}), 
\end{equation}
where $\mathcal{K} = e^{-\frac{c}{\ep}}\mu_1 \otimes \mu_2$ is the Gibbs kernel.
Since the entropic optimal transport~\eqref{KL reg OT} is equivalent to \eqref{Schrodinger bridge pb}, it is often referred to as the Schr\"{o}dinger bridge problem (see the survey papers \cite{chen2021stochastic, Leonard2014, Nutz2022introduction} and references therein). 

The problem of convex relaxation through entropy for linear programming has been extensively studied in the literature \cite{bregman1967relaxation, cominetti1994asymptotic, Erlander1980, takatsu2021relaxation, wilson1969}. A notable contribution was made by Cuturi~\cite{cuturi2013sinkhorn}, who showed that entropic optimal transport \eqref{KL reg OT} could be efficiently computed using the Sinkhorn--Knopp algorithm \cite{sinkhorn1964relationship, SK1967}. This breakthrough has led to the widespread application of entropic optimal transport in various fields such as image processing, graphics, and machine learning (see, for example, \cite{arjovsky2017wasserstein, patrini2020sinkhorn, solomon2015convolutional}). For a comprehensive overview of this topic, we refer interested readers to the book by Peyre--Cuturi~\cite{PeyreCutri} and the references therein.

The convergence rate of (\ref{KL reg OT}) as $\ep \to 0$ has attracted considerable interest from both mathematical and applied perspectives and has been the subject of active study. In the discrete setting, the convergence rate problem is studied by Cominetti--San Mart\'in \cite{cominetti1994asymptotic}. They showed that the convergence rate is exponentially fast. A quantitative estimate of the convergence rate has been studied by Weed \cite{weed2018}. 

In the continuum setting, Carlier--Duval--Peyr\'e--Schmitzer \cite{carlier2017convergence} proved $\Gamma$-convergence of \eqref{KL reg OT} to the Monge--Kantorovich problem \eqref{MK pb} in the case of a quadratic cost function (i.e., $c(x,y) = |x-y|^2$). As a result, they proved that $\mathrm{OT_{\ep}}$ converges to $\mathrm{OT}$. The block approximation method, introduced in this work, is a standard technique for estimating convergence rates in regularized optimal transport. The paper by Carlier--Pegon--Tamanini~\cite{carlier2023convergence} established the following convergence rate for \eqref{KL reg OT}: 
\eq{\label{esti KL reg OT}
  \mathrm{OT}_{\ep} - \mathrm{OT}
  \leq D \ep \log\frac{1}{\ep} + O(\ep), 
}
where the constant $D$ is the minimum of the entropy dimension of the marginals. Their approach is based on the block approximation method, and they also showed that the leading order term $\ep \log 1/\ep$ in the estimate \eqref{esti KL reg OT} is sharp. Recently, Eckstein--Nutz \cite{ENpre22, EN22} proposed a new method for obtaining convergence rates based on quantization and shadow arguments. This approach allows deriving the same results as in \eqref{esti KL reg OT} and can be applied to a wide range of regularizers, making it a highly useful tool. The methods are of great importance in this paper. 

Thus, the regularized transport problem with KL divergence has a long history and has been widely studied. Recently, the study of other regularizers has also attracted attention from the point of view of computational optimal transport. In particular, the regularization using Tsallis relative entropy, which is a one-parameter extension of KL divergence, has been studied in terms of sparsity and computational complexity of the entropic transport plan \cite{Bao2022,muzellec2017}. 

However, to our knowledge, few results focus on the properties of Tsallis regularization in the continuous setting and its differences with KL regularization, except in the case of multivariate normal distribution \cite{tong2021entropy}. In this paper, we consider the following regularized optimal transport problem with Tsallis relative entropy: 
\eq{\label{Tsallis reg OT}
  \mathrm{OT}_{q, \ep} = \inf_{\pi \in \Pi(\mu_1, \mu_2)} \left\{\int_{X_1 \times X_2} c\, d\pi + \ep D_{q}(\pi, \mu_1 \otimes \mu_2)\right\},  
}
where $D_q$ denotes the Tsallis relative entropy of the transport plan $\pi$ with respect to the product measure $\mu_1 \otimes \mu_2$, which is defined for general probability measures as 
\eq{\label{$D_q$-divergence}
  D_q(\mu, \nu) = \left\{
  \begin{aligned}
  &\frac{1}{q - 1} \int \left[\left(\frac{d\mu}{d\nu}\right)^q - \frac{d\mu}{d\nu}\right]\, d\nu
    && \mathrm{if} \ \mu \ll \nu \ \mathrm{and}\ q > 1, \\
  &\infty 
    && \mathrm{otherwise}. 
  \end{aligned}
  \right.
}
The relative entropy \eqref{$D_q$-divergence} is a generalization of the KL divergence and appears in so-called Tsallis statistical mechanics introduced by Tsallis \cite{T88}. It has been widely studied not only in nonlinear dynamics and nonequilibrium physics but also in information theory (see \cite{Borland1998, Furuichi2004, Furuichi2005, Naudts2011, Robledo2022, S04, tsallisbook2009} and references therein). 

We introduce $q$-logarithmic and $q$-exponential functions, which have been well studied in the field of Tsallis statistical mechanics. For a parameter $q \in \R$ with $q \neq 1$, the $q$-logarithmic and $q$-exponential functions are defined as 
\eqn{
  \log_q{y} = \frac{y^{1 - q} - 1}{1 - q}\quad (y > 0), \quad \exp_q{y} = [1 + (1-q)y]_+^{\frac{1}{1-q}}\quad (y \in \R),
}
where we define $f_+(x) = \max \left\{f(x), 0 \right\}$. 
They recover the natural definition of the logarithmic and exponential as $q \to 1$. One can represent $D_q$ by using the $q$-logarithmic function as follows:
\eqn{
  D_q(\mu, \nu)
  = \int \frac{d\mu}{d\nu} \log_{2 - q}\left(\frac{d\mu}{d\nu}\right)\, d\nu
}
for $\mu \ll \nu$ and $q > 1$. The purpose of this paper is to investigate the fundamental properties of Tsallis entropic regularized optimal transport \eqref{Tsallis reg OT} using the methods introduced by Eckstein--Nutz \cite{ENpre22, EN22}. Specifically, we employ their techniques to prove the $\Gm$-convergence as the regularization parameter $\ep$ approaches zero and demonstrate the convergence of the minimizer. Additionally, we establish the convergence rate of Tsallis entropic regularized optimal transport and compare it with the convergence rate of the KL divergence. By focusing on these aspects, we aim to establish a theoretical foundation for the application of Tsallis relative entropy in the context of optimal transport.

Let $p \in [1,\infty)$ and $\mathcal{P}_p(X)$ denote the set of probability measures over a Polish space $X$ with finite $p$-moment. Let $\mathcal{P}^n(X) \subset \mathcal{P}(X)$ denote the set of probability measures supported on at most $n$ points. We also assume that the cost function $c$ satisfies $(\mathrm{A}_{L, C})$ condition for some constants $L, C \geq 0$ introduced in \cite[Theorem 3.2]{ENpre22}, i.e., for any $\mu_i$, $\tilde{\mu}_i \in \mathcal{P}_p(X_i)$ ($i = 1, 2$) with $W_p(\tilde{\mu}_i, \mu_i) \leq C$, it holds that
\eq{\label{ineq;cost}
  \left|\int_{X_1 \times X_2} c\, d(\pi - \tilde{\pi})\right|
  \leq L W_p(\pi, \tilde{\pi})
}
for all $\pi \in \Pi(\mu_1, \mu_2)$ and $\tilde{\pi} \in \Pi(\tilde{\mu}_1, \tilde{\mu}_2)$, where $W_p(\pi, \tilde{\pi})$ is the $p$-th Wasserstein distance between $\pi$ and $\tilde{\pi}$. For instance, if the cost function $c$ is a $L$-Lipschitz function or  $c(x, y) = |x - y|^p$ on $X_1 \times X_2$ for $1 \leq p < \infty$, then $c$ satisfies $(\mathrm{A}_{L, C})$ condition. 
In the latter case, the constant $L$ depends on the moments of $\mu_i$ and on $C$. 

In the following, we present the main results of the paper. First, we consider the $\Gamma$-convergence of the Tsallis entropic OT \eqref{Tsallis reg OT}. Let us recall the definition of $\Gamma$-convergence. 

\begin{dfn}[$\Gamma$-convergence]\label{Gamma-convergence}
A sequence $\{F_k\}_{k \in \Nt}$ of functionals $F_k: \mathcal{P}_p(X_1 \times X_2) \to \R \cup \{ \infty \}$ is defined as 
\eqn{
  F_k(\pi) = \begin{cases} 
  \displaystyle \int_{X_1 \times X_2} c\, d\pi + \ep_k D_q(\pi, \mu_1 \otimes \mu_2) \quad &\text{if} \,\, \pi \in \Pi(\mu_1, \mu_2) \\
  +\infty &\text{otherwise}, 
  \end{cases}
}
where $\ep_k \to 0$ as $k \to \infty$, and also $F: \mathcal{P}_p(X_1 \times X_2) \to \R \cup \{ \infty \}$ defined as 
\eqn{
  F(\pi) = \begin{cases} 
  \displaystyle \int_{X_1 \times X_2} c\, d\pi \quad &\text{if} \,\, \pi \in \Pi(\mu_1, \mu_2) \\
  +\infty &\text{otherwise}. 
  \end{cases}
}

We say that the sequence $\{F_k\}_{k \in \Nt}$ $\Gm$-converges to $F$ with respect to the narrow topology if the following two conditions hold:
\begin{description}
\item[(liminf inequality)] For any sequence $\{\pi_k\}_{k \in \Nt} \subset \mathcal{P}_p(X_1 \times X_2)$ such that $\pi_k \to \pi$ narrowly, 
\eq{\label{ineq;liminf}
  F(\pi) \leq \liminf_{k \to \infty} F_k(\pi_k).
}
\item[(limsup inequality)] For any $\pi \in \mathcal{P}_p(X_1 \times X_2)$, there exists a sequence $\{\pi_k\}_{k \in \Nt} \subset \mathcal{P}_p(X_1 \times X_2)$ converging to $\pi_k \to \pi$ narrowly such that  
\eq{\label{ineq;limsup}
  F(\pi) \geq \limsup_{k \to \infty} F_k(\pi_k).
}
\end{description}
\end{dfn}

From the direct consequence of the definition of $\Gamma$-convergence, we obtain the following theorem. 

\begin{thm}\label{thm;narrow}
Let $\mu_i \in \mathcal{P}_p(X_i)$ ($i = 1, 2$) and $1 < q \leq 2$. Then, the sequence $\{F_k\}_{k \in \Nt}$ $\Gm$-converges to $F$ with respect to the narrow topology. Moreover, the unique minimizer $\pi^*_{k} \in \Pi(\mu_1, \mu_2)$ of $F_{k}$ converges to a minimizer $\pi^* \in \Pi(\mu_1, \mu_2)$ of $F$ narrowly. 
\end{thm}

Furthermore, by using the quantization and shadow arguments developed by Eckstein--Nutz \cite{EN22, ENpre22}, one can obtain the convergence rate of the Tsallis entropic regularized optimal transport \eqref{Tsallis reg OT}. To state the theorem, we recall the quantization for a probability measure, which comes from signal processing (see \cite{GL07}). Quantization for a probability measure is a topic related to approximating a $d$-dimensional probability measure $\mu$ with a given number $n$ of discrete probabilities in an optimal sense. 

\begin{dfn}[Quantization]
For $\mu \in \mathcal{P}_p(X)$, we say that $\mu$ satisfies $\mathrm{quant}_p(C, \bt)$ for some constants $C \geq 0$ and $\bt > 0$ if there exists $\{\mu^n\}_{n \in \Nt} \subset \mathcal{P}^n(X)$ such that
\eq{\label{ineq;quant}
  W_p(\mu^n, \mu)
  \leq C n^{- \frac{1}{\bt}}.
}
\end{dfn}

\begin{rem}
Let $X = \R^d$ and $\mu \in \mathcal{P}_{p + \dl}(\R^d)$ for some $\dl > 0$. Then one can take $\bt = d$ and $\mu$ satisfies $\mathrm{quant}_p(C, d)$ (see \cite[Theorem 6.2]{GL07}).  
\end{rem}

Then, the following theorem holds: 
\begin{thm}\label{thm;rate}
Let $\mu_i \in \mathcal{P}_p(X_i)$ ($i = 1, 2$). Suppose $\mu_2$ and the cost function $c$ satisfy $\mathrm{quant}_p(C, \bt)$ and $(\mathrm{A}_{L, C})$ condition for some constants $\bt \geq 1$ and $L, C \geq 0$, respectively. Then for $1 < q \leq 2$, 
\eq{\label{rate;Tsallis}
  \mathrm{OT}_{q, \ep} - \mathrm{OT}
  \leq K_{1, q} 
    \frac{\bt}{(q - 1) \bt + 1} \ep \log_\frac{1}{(q - 1) \bt + 1}{\frac{1}{\ep}}
    + \frac{K_{1, q} - 1}{q - 1} \ep
    + K_{2, q} \ep^\frac{1}{(q - 1) \bt + 1},
}
where
\eqn{
  K_{1, q} = \left(\frac{2 L C}{\bt}\right)^\frac{(q - 1) \bt}{(q - 1) \bt + 1}\quad\mathrm{and}\quad
  K_{2, q} = (2 L C)^\frac{(q - 1) \bt}{(q - 1) \bt + 1} \bt^\frac{1}{(q - 1) \bt + 1}.
}
\end{thm}

\begin{rem}\label{rem;rate}
We remark that from the result of \cite[Theorem 3.3]{ENpre22}, 
\eq{\label{rate;KL}
\mathrm{OT}_\ep - \mathrm{OT} \leq \bt \ep \log\frac{1}{\ep} + 4 L C \ep.
}
Since $q > 1$, we note that $\tilde{q} := \frac{1}{(q - 1) \bt + 1} < 1$. Thus we obtain that
\eqn{
  \lim_{x \to \infty} \frac{\log{x}}{\log_{\tilde{q}}{x}}
  = \lim_{x \to \infty} \frac{(1 - \tilde{q}) \log{x}}{x^{1 - \tilde{q}} - 1}
  = \lim_{x \to \infty} \frac{\frac{1}{x}}{x^{- \tilde{q}}}
  = \lim_{x \to \infty} x^{\tilde{q} - 1} 
  = 0. 
}
Therefore, when considering $\ep \to 0$, the convergence rate \eqref{rate;KL} is faster than \eqref{rate;Tsallis}. Since the KL divergence \eqref{KL divergence} is the limit case of Tsallis relative entropy \eqref{$D_q$-divergence}, we can conclude that KL is the fastest convergence rate in terms of Tsallis relative entropy.  
\end{rem}

\begin{rem}
We compute the explicit constants in \eqref{rate;Tsallis} while similar results to \eqref{rate;Tsallis} have been obtained in \cite{ENpre22}, where they considered the $L^\rho$ regularization of the optimal transport for $\rho > 1$ (see Example 3.1 in \cite{ENpre22}). 
We remark that the estimates from $L^\rho$ regularization do not coincide with \eqref{rate;KL} as $\rho \to 1$, while in our setting, \eqref{rate;KL} can be derived from \eqref{rate;Tsallis} when $q \to 1$. In fact, we observe that $K_{1, q} \to 1$ and $K_{2, q} \to \bt$ as $q \to 1$. Furthermore, the convergence of the logarithmic function of the $q$-logarithmic function as $q \to 1$ implies that the leading term of the right-hand side of the convergence rate of the Tsallis entropic regularization \eqref{rate;Tsallis} converges to $\bt \ep \log (\frac{1}{\ep})$. Therefore, when $q \to 1$, the leading term of the Tsallis entropic regularization \eqref{rate;Tsallis} is of the same order as the leading term of the KL divergence \eqref{rate;KL}.
\end{rem}

Moreover, one can prove the sharpness of the convergence rate of the Tsallis entropic regularized optimal transport \eqref{rate;Tsallis}. 
\begin{thm}[Sharpness of convergence rate]\label{thm;sharp}
Let $X_1 = X_2 = \R^d$ with $\mu_1 = \mu_2$ the uniform distribution on $[0, 1]^d$ and $c(x, y) = \sum_{i = 1}^d |x_i - y_i|$.  
Then for $q > 1$, 
\eqn{
  \mathrm{OT}_{q, \ep} - \mathrm{OT}
  \geq \tilde{K}_{1, q} \ep^\frac{1}{(q - 1) d + 1} - \tilde{K}_{2, q} \ep, 
}
where
\begin{align*}
  \tilde{K}_{1, q}
  &= \frac{(q - 1) d + 1}{(q - 1) d + q} \left[\frac{q - 1}{((q - 1) d + q) 2^\frac{(q - 1) d + 1}{q - 1}} \left(\frac{q}{q - 1}\right)^\frac{q}{q - 1}\right]^\frac{q - 1}{(q - 1) d + 1}, \\ 
  \tilde{K}_{2, q} &= \left(1 + 2^\frac{(q - 1) d + 1}{q - 1}\right) q^{- \frac{q}{q - 1}}.
\end{align*}
\end{thm}

The paper is organized as follows. In Section \ref{preliminaries}, we recall the data processing inequality and the definition and properties of the shadow of a coupling introduced by Eckstein--Nutz \cite{EN22, ENpre22}. In Section \ref{pf of Gamma conv}, we prove $\Gamma$-convergence of the Tsallis entropic regularized optimal transport \eqref{Tsallis reg OT} by using the quantization and shadow arguments. In Section \ref{convergence rate of Tsallis}, we show the convergence rate of the Tsallis entropic regularized optimal transport \eqref{rate;Tsallis}. One can see that the convergence rate of the optimal transportation cost with KL regularization \eqref{rate;KL} converges faster than with Tsallis entropic regularization \eqref{rate;Tsallis}. In Section \ref{sharpness of convergence rate}, we prove the sharpness of the convergence rate obtained in Section \ref{convergence rate of Tsallis}. 

\section{Preliminaries}\label{preliminaries}
In what follows, we introduce the data processing inequality and explain the argument of the shadow of a coupling introduced by Eckstein--Nutz \cite{EN22, ENpre22}. 

For simplicity, we write
\eqn{
  f_q(x) 
  = \left\{
    \begin{aligned}
    &\frac{x^q - x}{q - 1} && \mathrm{if}\ q > 1, \\
    &x \log{x} && \mathrm{if}\ q = 1
    \end{aligned}
    \right. \quad \mathrm{and}\quad
  \phi_q(x) 
  = \frac{f_q(x)}{x} = \left\{
    \begin{aligned}
    &\frac{x^{q - 1} - 1}{q - 1} && \mathrm{if}\ q > 1, \\
    &\log{x} && \mathrm{if}\ q = 1.
    \end{aligned}
    \right.
}
Then one can represent the Tsallis relative entropy \eqref{$D_q$-divergence} by
\eqn{
  D_q(\mu, \nu)
  = \int f_q\left(\frac{d\mu}{d\nu}\right)\, d\nu
  = \int \frac{d\mu}{d\nu}
    \phi_q\left(\frac{d\mu}{d\nu}\right)\, d\nu
    = \int \phi_q\left(\frac{d\mu}{d\nu}\right)\, d\mu. 
}
We note that $\phi_q$ is concave when $1 \leq q \leq 2$, while $f_q$ is convex for any $q \geq 1$.

\subsection{Markov kernels and the data processing inequality}
We explain the Markov kernels (also called the stochastic kernels) and the data processing inequality. These concepts and inequality are important tools in information theory and play an important role in controlling the Tsallis relative entropy \eqref{$D_q$-divergence} in this paper. We refer the interested reader to \cite{cinlar2011probability, polyanskiy2022information}. 
\begin{dfn}[Markov kernel]
Let $X_1$ and $X_2$ be Polish spaces. The mapping $K$ is called the Markov kernel if the following two conditions hold:  
\begin{enumerate}
    \item For every measurable set $B$, the mapping $x \mapsto K(x,B)$ is a measurable function on $X_1$.
    \item For every $x \in X_1$, the mapping $B \mapsto K(x,B)$ is a probability measure on $X_2$. 
\end{enumerate}
\end{dfn}
Markov kernels can be interpreted as random operations or transformations. 
A probability measure $\pi \in \mathcal{P}(X_1 \times X_2)$ can be disintegrated by a Markov kernel. According to the result of \cite[Chapter IV, Theorem 2.18]{cinlar2011probability}, there exists a probability measure $\mu_1 \in \mathcal{P}(X_1)$ and a Markov kernel $K: X_1 \to \mathcal{P}(X_2)$ such that 
\begin{equation}\label{disinteg}
    \pi(A \times B) = \int_A \mu_1(dx) K(x,B) \quad A \subset X_1, \,\, B \subset X_2. 
\end{equation}

We denote the second marginal of $\pi = \mu_1 \otimes K\in \mathcal{P}(X_1 \times X_2)$ by $\mu_1 K \in \mathcal{P}(X_2)$. 
Utilizing the disintegration of the coupling, it is possible to obtain the bounds of the Tsallis relative entropy \eqref{$D_q$-divergence}. 
\begin{lem}\label{lem;diverg}
Let $\pi \in \Pi(\mu_1, \mu_2)$. If $\mu_2 \in \mathcal{P}^n(X_2)$, then for $1 < q \leq 2$,
\eq{\label{ineq;diverg}
  D_q(\pi, \mu_1 \otimes \mu_2)
  \leq \phi_q(n).
}
\end{lem}

{\it Proof of Lemma \ref{lem;diverg}}
Let $\s$ be the counting measure of the support of $\mu_2$ and $\{x_i\}_{i  = 1}^n$ be the support of $\mu_2$. Since $\mu_2 \in \mathcal{P}^n(X_2)$, we may assume that $\s = \Sigma^n_{i=1} \delta_{x_i}$ and $\mu_2 = \Sigma^n_{i=1} \alpha_i \delta_{x_i}$, where $\Sigma^n_{i=1} \alpha_i = 1$. Consider the disintegration of $\pi$ as $\mu_1 \otimes K$. Since $K: X_1 \to \mathcal{P}^n(X_2)$, we can write $K = \Sigma^n_{i=1} \beta_i \delta_{x_i}$, where $\Sigma^n_{i=1} \beta_i = 1$. 
Then we have 
\eq{\label{ineq;temp1}
  \frac{d\pi}{d(\mu_1 \otimes \mu_2)}
  = \frac{d\mu_1}{d\mu_1}\frac{dK}{d\mu_2}
  = \frac{dK}{d\mu_2}
  = \frac{\beta_i}{\alpha_i}
  \leq \frac{1}{\alpha_i} = \frac{d\s}{d\mu_2}.
}
Since the function $\phi_q$ is non-decreasing, by \eqref{ineq;temp1} we have   
\eqn{
  D_q(\pi, \mu_1 \otimes \mu_2)
  = \int_{X_1 \times X_2} \phi_q\left(\frac{d\pi}{d(\mu_1 \otimes \mu_2)}\right)\, d\pi
  \leq \int_{X_1 \times X_2} \phi_q\left(\frac{d\s}{d\mu_2}\right)\, d\pi.
}
The Jensen inequality yields from the concavity of $\phi_q$ that 
\begin{align*}
  D_q(\pi, \mu_1 \otimes \mu_2)
  &\leq \int_{X_1 \times X_2} \phi_q\left(\frac{d\s}{d\mu_2}\right)\, d\pi
  = \int_{X_2} \phi_q\left(\frac{d\s}{d\mu_2}\right)\, d\mu_2 \\
  &\leq \phi_q\left(\int_{X_2}\, d\s\right) = \phi_q(n),
\end{align*}
which conclude that the inequality \eqref{ineq;diverg} holds.

\begin{lem}[Data processing inequality]\label{lem;data}
Let $\mu, \nu \in \mathcal{P}(X_1)$ and let $K: X_1 \to \mathcal{P}(X_2)$ be a Markov kernel. Then the following inequality for the Tsallis relative entropy \eqref{$D_q$-divergence} holds: 
\eq{\label{ineq;data}
  D_q(\mu K, \nu K)
  \leq D_q(\mu, \nu).
}
\end{lem}

{\it Proof of Lemma \ref{lem;data}}
This proof follows the same manner in {\cite[Lemma 4.1]{EN22}} for the Tsallis relative entropy \eqref{$D_q$-divergence}. We may assume that $\mu \ll \nu$. For any Markov kernels $K_1 \ll K_2: X_1 \to \mathcal{P}(X_2)$, we have
\eqn{
  \frac{d (\mu \otimes K_1)}{d (\nu \otimes K_2)}
  = \frac{d\mu}{d\nu}\frac{dK_1}{dK_2}
    \quad \nu \otimes K_2\text{-a.s.}
}
We note that
\eq{\label{data;1}
\begin{split}
  D_q(\mu \otimes K, \nu \otimes K)
  &= \int_{X_1 \times X_2} f_q\left(\frac{d(\mu \otimes K)}{d(\nu \otimes K)}\right) d (\nu \otimes K) \\
  &= \int_{X_1} f_q\left(\frac{d\mu}{d\nu}\right) d\nu
  = D_q(\mu, \nu),  
\end{split}
}
where we used $\displaystyle \int_{X_2} dK = 1$. Since $f_q$ is a convex function if $q > 1$, then by the Jensen inequality, we have
\begin{align}
  D_q(\mu \otimes K_1, \nu \otimes K_2)
  &= \int_{X_1 \times X_2} f_q\left(\frac{d\mu}{d\nu} \frac{dK_1}{dK_2}\right) d(\nu \otimes K_2) \notag \\
  &\geq \int_{X_1} f_q\left(\frac{d\mu}{d\nu}\right)\, d \nu
  = D_q(\mu, \nu). \label{data;2}
\end{align}
Consider the disintegrations of $\mu \otimes K$ and $\nu \otimes K$ as $(\mu K) \otimes \tilde{K}_1$ and $(\nu K) \otimes \tilde{K}_2$ from the second marginal to the first, respectively. Applying \eqref{data;2}, we have
\eq{\label{data;3}
  D_q(\mu \otimes K, \nu \otimes K)
  = D_q\big((\mu K) \otimes \tilde{K}_1, (\nu K)\otimes \tilde{K}_2\big)
  \geq D_q(\mu K, \nu K)
}
Combining these inequalities \eqref{data;1} and \eqref{data;2}, we obtain the inequality \eqref{ineq;data}.
\qed

\subsection{Shadows}
The idea of the shadow is to consider a proper $W_p$ projection of a coupling $\pi$ onto the set of couplings of other marginals.

\begin{dfn}[{\cite[Definition 3.1]{EN22}}]\label{dfn;shadow}
Let $1 \leq p < \infty$ and $\mu_i, \tilde{\mu}_i \in \mathcal{P}_p(X_i)$ ($i = 1, 2$).  
Let $\kp_i \in \Pi(\mu_i, \tilde{\mu}_i)$ be a optimal coupling for $W_p(\mu_i, \tilde{\mu}_i)$, and let $\kp_i = \mu_i \otimes K_i$ be a disintegration.  
For $\pi \in \Pi(\mu_1, \mu_2)$, its shadow $\tilde{\pi} \in \Pi(\tilde{\mu}_1, \tilde{\mu}_2)$ is defined as the second marginal of $\pi \otimes K \in \mathcal{P}_p((X_1 \times X_2)^2)$, where the kernel $K: X_1 \times X_2 \to \mathcal{P}_p(X_1 \times X_2)$ is defined as $K(x) = K_1(x_1) \otimes K_2(x_2)$.  
\end{dfn}

Recall the properties of the shadow. By using the argument of the shadow and the data processing inequality \eqref{ineq;data}, one can control the Tsallis relative entropy $D_q$. These properties play an important role throughout this paper. 

\begin{lem}[{\cite[Lemma 3.2]{EN22}}]\label{lem;shadow}
Let $1 \leq p < \infty$ and $\mu_i, \tilde{\mu}_i \in \mathcal{P}_p(X_i)$ ($i = 1, 2$).  
For any $\pi \in \Pi(\mu_1, \mu_2)$, 
its shadow $\tilde{\pi} \in \Pi(\tilde{\mu}_1, \tilde{\mu}_2)$ satisfies the following:
\eq{
  W_p(\pi, \tilde{\pi})^p
  = W_p(\mu_1, \tilde{\mu}_1)^p + W_p(\mu_2, \tilde{\mu}_2)^p,
}
\eq{
  D_q(\tilde{\pi}, \tilde{\mu}_1 \otimes \tilde{\mu}_2)
  \leq D_q(\pi, \mu_1 \otimes \mu_2).
}
\end{lem}

\section{$\Gamma$-convergence of the Tsallis entropic regularized optimal transport}\label{pf of Gamma conv}

In this section, we prove the $\Gamma$-convergence of the Tsallis entropic regularized optimal transport \eqref{Tsallis reg OT}. 
The proof of Theorem \ref{thm;narrow} can be divided into the following two propositions (liminf inequality \eqref{ineq;liminf} and limsup inequality \eqref{ineq;limsup}): 

\begin{prop}\label{prop;liminf}
For any sequence $\{\pi_k\}_{k \in \Nt} \subset \mathcal{P}_p(X_1 \times X_2)$ such that $\pi_k \to \pi$ narrowly, we have
\eqn{
  F(\pi) \leq \liminf_{k \to \infty} F_k(\pi_k).
}
\end{prop}

{\it Proof of Proposition \ref{prop;liminf}}
The Jensen inequality implies that 
\begin{align*}
D_q(\pi_k, \mu_1 \otimes \mu_2) &= \int_{X_1 \times X_2} f_q\left(\frac{d\pi_k}{d(\mu_1 \otimes \mu_2)}\right) \, d(\mu_1 \otimes \mu_2) \\
&\geq f_q \left(\int_{X_1 \times X_2} d\pi_k \right) = f_q(1) = 0. 
\end{align*}
Therefore, we obtain
\eq{\label{ineq_liminf;divergence}
  \liminf_{k \to \infty} \ep_k D_q(\pi_k, \mu_1 \otimes \mu_2)
  \geq 0.
}
Moreover, it follows from the lower semi-continuity of the coupling (see \cite[Lemma 4.3]{villani2009optimal}) that
\eq{\label{ineq_liminf;lower semi}
  \int_{X_1 \times X_2} c\, d\pi
  \leq \liminf_{k \to \infty} \int_{X_1 \times X_2} c\, d\pi_k.
}
Combining \eqref{ineq_liminf;divergence} with \eqref{ineq_liminf;lower semi}, the liminf inequality \eqref{ineq;liminf} holds.  
\qed

\begin{prop}\label{prop;limsup}
Let $\mu_i \in \mathcal{P}_p(X_i)$ ($i = 1, 2$) and $\{\ep_k\}_{k \in \Nt}$ be a sequence of positive numbers converging zero and satisfying $\ep_k \phi_q(k) \to 0$. Assume that $1 < q \leq 2$. For any $\pi \in \mathcal{P}_p(X_1 \times X_2)$, there exists a sequence $\{\pi_k\}_{k \in \Nt} \subset \mathcal{P}_p(X_1 \times X_2)$ converging to $\pi_k \to \pi$ narrowly such that  
\eqn{
  F(\pi) \geq \limsup_{k \to \infty} F_k(\pi_k).
}
\end{prop}

{\it Proof Proposition \ref{prop;limsup}}
By the assumption and the argument of quantization, there exists $\{ \mu^k_2 \}_{k \in \Nt} \subset \mathcal{P}_p^k$ such that 
\eqn{
  \lim_{k \to \infty} W_p(\mu^k_2, \mu_2) = 0
}
(see \cite[Theorm 6.18]{villani2009optimal}).
For any $\pi \in \Pi(\mu_1, \mu_2)$, we take the shadow $\tilde{\pi} \in \Pi(\mu_1, \mu_2^k)$ of $\pi$. 
Taking the shadow $\pi_k \in \Pi(\mu_1, \mu_2)$ of $\tilde{\pi}$ again, it follows from the triangle inequality concerning the Wasserstein distance  
(see \cite[Theorem 7.3]{villani2021topics}) and Lemma \ref{lem;shadow} that
\eqn{
  W_p(\pi_k, \pi)
  \leq W_p(\pi_k, \tilde{\pi}) + W_p(\tilde{\pi}, \pi)
  \leq 2 W_p(\mu_2^k, \mu_2)
  \to 0 \,\, \mathrm{as} \,\, k \to \infty.
}
Thus, one can see from a property of the Wasserstein distance (see \cite[Theorem 7.12]{villani2021topics}) that
\eqn{
  \pi_k \to \pi \,\, \text{narrowly as} \,\, k \to \infty.
}
Also, since the cost function $c$ satisfies $(\mathrm{A}_{L, C})$ condition, we have 
\eqn{
  \left|\int_{X_1 \times X_2} c\, d (\pi_k - \pi)\right| \leq L W_p(\pi_k, \pi) \to 0 \,\, \mathrm{as} \,\, k \to \infty.
}
Therefore, we obtain 
\eq{\label{conv cost}
  \int_{X_1 \times X_2} c\, d\pi_k \to \int_{X_1 \times X_2} c\, d\pi \,\, \mathrm{as}\,\, k \to \infty.
}
On the other hand, Lemma \ref{lem;diverg} and Lemma \ref{lem;shadow} imply that 
\eqn{
  D_q(\pi_k, \mu_1 \otimes \mu_2) \leq D_q(\tilde{\pi}, \mu_1 \otimes \mu^k_2) \leq \phi_q(k).  
}
Then, we obtain 
\eq{\label{conv tsallis}
  \limsup_{k \to \infty} \left( \ep_k D_q(\pi_k, \mu_1 \otimes \mu_2) \right)
  \leq \limsup_{k \to \infty} \left(\ep_k \phi_q(k) \right) = 0. 
}
By combining \eqref{conv cost} and \eqref{conv tsallis}, limsup inequality \eqref{ineq;limsup} holds. 
\qed

{\it Proof of Theorem \ref{thm;narrow}}
Proposition \ref{prop;liminf}, \ref{prop;limsup} imply that $\{F_k\}_{k \in \Nt}$ $\Gm$-converges to $F$ with respect to narrow topology. Moreover, by the strict convexity of Tsallis relative entropy \eqref{$D_q$-divergence}, $\pi^{*}_k$ is the unique minimizer of $F_k$. Since $\Pi(\mu_1,\mu_2)$ is tight (see \cite[Lemma 4.4]{villani2009optimal}), the sequence $\{\pi^{*}_k\}_{k \in \Nt} \subset \Pi(\mu_1,\mu_2)$ converges to a coupling $\pi^{*} \in \Pi(\mu_1,\mu_2)$ narrowly. By using the property of $\Gm$-convergence, $\pi^{*}$ is a minimizer of $F$. This is the desired conclusion. 
\qed

\section{Convergence rate of the Tsallis entropic regularized optimal transport}\label{convergence rate of Tsallis}

In this section, we prove Theorem \ref{thm;rate}. In the proof of Theorem \ref{thm;rate}, we use the double shadow argument of the coupling developed by the paper \cite{ENpre22,EN22}. 

{\it Proof of Theorem \ref{thm;rate}}
Let $\pi^* \in \Pi(\mu_1, \mu_2)$ be the optimal coupling of  $\mathrm{OT}$. Since $\mu_2$ satisfies $\mathrm{quant}_p(C, \bt)$, there exists $\{\mu_2^n\}_{n \in \Nt} \subset \mathcal{P}^n(X_2)$ such that
\eq{\label{rate;1}
  W_p(\mu_2^n, \mu_2)
  \leq C n^{- \frac{1}{\bt}}.
}

For $n \in \Nt$, we take a coupling $\pi \in \Pi(\mu_1, \mu_2)$ of $\pi^*$ as a double shadow: 
Let $\tilde{\pi} \in \Pi(\mu_1, \mu_2^n)$ be a shadow of $\pi^*$ and define $\pi \in \Pi(\mu_1, \mu_2)$ as a shadow of $\tilde{\pi}$.  
By the triangle inequality and Lemma \ref{lem;shadow},
we have
\eq{\label{rate;2}
  W_p(\pi, \pi^*)
  \leq W_p(\pi, \tilde{\pi}) + W_p(\tilde{\pi}, \pi^*)
  \leq 2 W_p(\mu_2^n, \mu_2).
}
Since the cost function $c$ satisfies $(\mathrm{A}_{L, C})$ condition, it follows from \eqref{rate;1} and \eqref{rate;2} that
\eq{\label{rate;3}
  \int_{X_1 \times X_2} c\, d\pi - \int_{X_1 \times X_2} c\, d\pi^*
  \leq 2 L W_p(\mu_2^n, \mu_2)
  \leq 2 L C n^{- \frac{1}{\bt}}.
}
On the other hand, since $\pi$ is a shadow of $\tilde{\pi}$, Lemma \ref{lem;shadow} again yields
\eq{\label{rate;4}
  D_q(\pi, \mu_1 \otimes \mu_2)
  \leq D_q(\tilde{\pi}, \mu_1 \otimes \mu_2^n).
}
Thus, combining \eqref{rate;3} and \eqref{rate;4}, we have
\eqn{
  \mathrm{OT}_{q, \ep} - \mathrm{OT}
  \leq 2 L C n^{- \frac{1}{\bt}} + \ep D_q(\tilde{\pi}, \mu_1 \otimes \mu_2^n).
}
By Lemma \ref{lem;diverg}, we see that
\eqn{
  D_q(\tilde{\pi}, \mu_1 \otimes \mu_2^n)
  \leq \phi_q(n),
}
and hence, we obtain
\eq{\label{rate;5}
  \mathrm{OT}_{q, \ep} - \mathrm{OT}
  \leq 2 L C n^{- \frac{1}{\bt}} + \ep \phi_q(n).
}

Then the inequality \eqref{rate;5} implies that 
\eqn{
  \mathrm{OT}_{q, \ep} - \mathrm{OT}
  \leq 2 L C n^{- \frac{1}{\bt}} + \frac{\ep}{q - 1} (n^{q - 1} - 1).
}
If we take the natural number $n \in \Nt$ as
\eqn{
  n \sim \left(\frac{2 L C}{\bt}\right)^\frac{\bt}{(q - 1) \bt + 1} \ep^{- \frac{\bt}{(q - 1) \bt + 1}},
}
then we obtain the convergence rate of the Tsallis entropic regularized optimal transport \eqref{Tsallis reg OT}: 
\eq{\label{rate;Tsallis_temp}
\mathrm{OT}_{q, \ep} - \mathrm{OT} \leq \frac{1}{q - 1} 
    \left(K_{1, q} \ep^\frac{1}{(q - 1) \bt + 1} - \ep\right) + K_{2, q} \ep^\frac{1}{(q - 1) \bt + 1}, 
}
where the constants $K_{1, q}$ and $K_{2, q}$ are defined as 
\eqn{
  K_{1, q} = \left(\frac{2 L C}{\bt}\right)^\frac{(q - 1) \bt}{(q - 1) \bt + 1}\quad\mathrm{and}\quad
  K_{2, q} = (2 L C)^\frac{(q - 1) \bt}{(q - 1) \bt + 1} \bt^\frac{1}{(q - 1) \bt + 1}.
}
When we use the $q$-logarithmic function in the first term on the right-hand side of the convergence rate \eqref{rate;Tsallis_temp}, one can rewrite as follows:
\begin{align*}
  &\frac{1}{q - 1} 
    \left(K_{1, q} \ep^\frac{1}{(q - 1) \bt + 1} - \ep\right) \\
  &\quad = K_{1, q} \frac{\ep}{q - 1} 
    \left(\ep^{- \left(1 - \frac{1}{(q - 1) \bt + 1}\right)} - 1\right)
    + \frac{\ep}{q - 1} (K_{1, q} - 1) \\
  &\quad = K_{1, q} 
    \frac{1 - \frac{1}{(q - 1) \bt + 1}}{q - 1} \frac{\ep}{1 - \frac{1}{(q - 1) \bt + 1}}
    \left(\ep^{- \left(1 - \frac{1}{(q - 1) \bt + 1}\right)} - 1\right)
    + \frac{\ep}{q - 1} (K_{1, q} - 1) \\
  &\quad = K_{1, q} 
    \frac{\bt}{(q - 1) \bt + 1} \ep \log_\frac{1}{(q - 1) \bt + 1}{\frac{1}{\ep}}
    + \frac{\ep}{q - 1} (K_{1, q} - 1).
\end{align*}
This is the desired conclusion. 
\qed

\section{Sharpness of the convergence rate of the Tsallis entropic regularized optimal transport}\label{sharpness of convergence rate}
In this section, we prove the sharpness of the convergence rate of the Tsallis entropic regularized optimal transport \eqref{rate;Tsallis}. 
We introduce the dual problem of the entropic optimal transport problem to consider the sharpness of the convergence rate (see \cite[Theorem 3.6]{di2020optimal}).

For a convex function $f$, the Legendre transformation of $f$ is defined by
\eqn{
  f^*(y) = \sup_{ x \geq 0} \left\{x y - f(x) \right\}, \,\, y \in \R.
}
By considering the Legendre transformation of $f_q$, one can represent the dual problem of $\mathrm{OT}_{q, \ep}$ as 
\eq{\label{dual}
  \begin{split}
  \mathrm{OT}_{q, \ep} 
  &= \sup_{\substack{(h_1, h_2) \in L^1(\mu_1) \times L^1(\mu_2)}} 
    \left\{\int_{X_1} h_1(x)\, d\mu_1(x) + \int_{X_2} h_2(y) \, d\mu_2(y) \right. \\
    &\qquad\qquad\qquad \left. 
    - \ep \int_{X_1 \times X_2} f_q^*\left(\frac{h_1(x) + h_2(y) - c(x, y)}{\ep}\right)\, d\mu_1(x)\, d\mu_2(y)\right\},
  \end{split}
}
where $f_q^*$ is given as 
\eqn{
  f_q^*(y)
  = \left(\frac{q - 1}{q} y + \frac{1}{q}\right)^\frac{q}{q - 1} 
  = q^{- \frac{q}{q - 1}} \left[ 1 + (q - 1) y \right]_+^\frac{q}{q - 1} 
  = q^{- \frac{q}{q - 1}} \exp_{2 - q}(y)^q
}
for $y \in \R$. Also, the dual problem of $\mathrm{OT}$ is 
\begin{align*}
\mathrm{OT} = \sup_{\substack{(h_1, h_2) \in L^1(\mu_1) \times L^1(\mu_2), \\ h_1 + h_2 \leq c}} \left\{\int_{X_1} h_1(x)\, d\mu_1(x) + \int_{X_2} h_2(y) \, d\mu_2(y) \right\}.
\end{align*}

{\it Proof of Theorem \ref{thm;sharp}}
Since $X_1 = X_2 = \R^d$ and $\mu_1 = \mu_2$ is the uniform distribution on $[0, 1]^d$, 
we see that $\mathrm{OT} = 0$.
If we set $h_1 = a$ and $h_2 = 0$ for any $a \in \R$ in the dual problem \eqref{dual}, 
then it holds that
\begin{align*}
  \mathrm{OT}_{q, \ep} - \mathrm{OT} 
  &\geq \int_{X_1} h_1(x)\, d\mu_1(x) + \int_{X_2} h_2(y) \, d\mu_2(y) \\
    &\quad - \ep \int_{X_1 \times X_2} f_q^*\left(\frac{h_1(x) + h_2(y) - c(x, y)}{\ep}\right)\, d\mu_1(x)\, d\mu_2(y) \\
  &= a - \ep \int_{X_1 \times X_2} f_q^*\left(\frac{a - c(x, y)}{\ep}\right)
    \, d\mu_1(x)\, d\mu_2(y). 
\end{align*}
Thus, we have
\eq{\label{sharp;1}
  \mathrm{OT}_{q, \ep} - \mathrm{OT}
  \geq \sup_{a \in \R} \left\{a - \ep \int_{X_1 \times X_2} f_q^*\left(\frac{a - c(x, y)}{\ep}\right)
    \, d\mu_1(x)\, d\mu_2(y)\right\}.
}
Since $f^*_q$ is not decreasing and $c(x, y) \geq 0$, we have
\begin{align}
  &\ep \int_{X_1 \times X_2} f_q^*\left(\frac{a - c(x, y)}{\ep}\right)\, d\mu_1(x)\, d\mu_2(y) \notag \\
  &= \ep \int_{X_1 \times X_2} f_q^*\left(\frac{a - c(x, y)}{\ep}\right) \chi_{\left\{a \geq c(x, y)\right\} }\, d\mu_1(x)\, d\mu_2(y) \notag \\
    &\qquad + \ep \int_{X_1 \times X_2} f_q^*\left(\frac{a - c(x, y)}{\ep}\right)
    \chi_{\left\{a < c(x, y)\right\}}\, d\mu_1(x)\, d\mu_2(y) \notag \\
  &\leq \ep f_q^*\left(\frac{a}{\ep}\right) \int_{X_1 \times X_2} \chi_{\left\{a \geq c(x, y)\right\}}
    \, d\mu_1(x)\, d\mu_2(y) + \ep f^*_q(0), \label{sharp;2}
\end{align}
where $\chi_A$ is the characteristic function of the set $A$. Therefore, by putting \eqref{sharp;2} to \eqref{sharp;1}, we obtain 
\eqn{
  \mathrm{OT}_{q, \ep} - \mathrm{OT}
  \geq \sup_{a \in \R} \left\{a - \ep f^*_q\left(\frac{a}{\ep}\right) \int \chi_{\left\{a \geq c(x, y)\right\}}
    \, d\mu_1(x)\, d\mu_2(y)\right\}
    - q^{- \frac{q}{q - 1}} \ep.
}
By the definition of the cost function $c$, one can see that
\eqn{
  \chi_{ \left\{a \geq c(x, y)\right\} }
  \leq \prod_{i = 1}^d \chi_{ \left\{|x_i - y_i| \leq a \right\} }.
}
Thus, it holds that
\begin{align}
  \int_{X_1 \times X_2} \chi_{ \left\{a \geq c(x,y) \right\}}\, d\mu_1(x) d\mu_2(y)
  &\leq \prod_{i = 1}^d \int_0^1 \int_0^1 \chi_{\{|x_i - y_i| \leq a\}}\, dx_i\, dy_i \notag \\
  &= (2 a - a^2)^d 
  \leq (2 a)^d \label{sharp;3}
\end{align}
for any $0 \leq a \leq 1$. On the other hand,  
\eq{\label{sharp;4}
  \ep f^*_q\left(\frac{a}{\ep}\right)
  \leq \ep 2^{\frac{q}{q - 1} - 1} \left[q^{- \frac{q}{q - 1}} + \left(\frac{a (q - 1)}{\ep q}\right)^\frac{q}{q - 1}\right]
  \leq C_1 \ep^{- \frac{1}{q - 1}} a^\frac{q}{q - 1} + C_2 \ep,
}
where 
\eqn{
  C_1 = 2^\frac{1}{q - 1} \left(\frac{q - 1}{q}\right)^\frac{q}{q - 1}\quad \mathrm{and}\quad 
  C_2 = 2^\frac{1}{q - 1} q^{- \frac{q}{q - 1}}.
}
Hence, we see from \eqref{sharp;3} and \eqref{sharp;4} that 
\begin{align*}
  &a - \ep f^*\left(\frac{a}{\ep}\right) \int_{X_1 \times X_2} \chi_{\left\{a \geq c(x, y)\right\}}
    \, d\mu_1(x)\, d\mu_2(y) \\
  &\quad \geq a - 2^d C_1 \ep^{- \frac{1}{q - 1}} a^\frac{(q - 1) d + q}{q - 1} - 2^d C_2 \ep a^d.
\end{align*}
Choosing 
\eqn{
  a = k \ep^\frac{1}{(q - 1) d + 1}\quad \mathrm{with}\  k = \left[\frac{q - 1}{((q - 1) d + q) 2^d C_1}\right]^\frac{q - 1}{(q - 1) d + 1},
}
then we obtain
\eqn{
  \mathrm{OT}_{q, \ep} - \mathrm{OT}
  \geq \tilde{K}_{1, q} \ep^\frac{1}{(q - 1) d + 1} - \tilde{K}_{2, q} \ep,
}
where 
\eqn{
  \tilde{K}_{1, q}
  = \frac{(q - 1) d + 1}{(q - 1) d + q} \left[\frac{q - 1}{((q - 1) d + q) 2^d C_1}\right]^\frac{q - 1}{(q - 1) d + 1} \ \mathrm{and}\
  \tilde{K}_{2, q} = q^{- \frac{q}{q - 1}} + 2^d C_2. 
}
This is the desired conclusion. 
\qed

\section{Conclusions}
In this paper, we studied the continuous Tsallis entropic regularized optimal transport problem \eqref{Tsallis reg OT}. We discussed the $\Gamma$-convergence to the Monge--Kantorovich problem \eqref{MK pb} as the regularization parameter $\varepsilon$ approaches zero (Theorem \ref{thm;narrow}) and established the convergence rate (Theorem \ref{thm;rate}). In particular, we compared the convergence rate with that of the KL entropic regularized optimal transport problem \eqref{KL reg OT} and proved that the rate is slower in the case $1 < q \leq 2$. Additionally, we demonstrated that this convergence rate is sharp (Theorem \ref{thm;sharp}).

This paper represents a first step in the theoretical analysis of the continuous Tsallis entropic regularized optimal transport problem \eqref{Tsallis reg OT}, leaving many open questions. For instance, for $0 < q < 1$, the convergence rate of the Tsallis entropic regularized optimal transport problem \eqref{Tsallis reg OT} is expected to be faster than that of the KL entropic regularized optimal transport problem \eqref{KL reg OT}, as discussed in Remark~\ref{rem;rate}. However, due to the technical dependence on the concavity of $\phi_q$, the convergence rate for $0 < q < 1$ remains unresolved. In the discrete case, the convexity of the relative entropy significantly influences the convergence rate in the regularized optimal transport problem in terms of the Bregman divergence \cite{morikuni2023error}. A similar situation is expected for the continuous case, which will be explored in future work.

Furthermore, it would be interesting to quantitatively estimate the spread of the support of the regularized transport plan in Tsallis relative entropy from the perspective of computational optimal transport. The case of KL entropy is discussed in \cite{bernton2022entropic}.

As a natural generalization of the Schr\"{o}dinger bridge problem \eqref{Schrodinger bridge pb} in the context of Tsallis relative entropy, we consider the following minimization problem:
\begin{equation}\label{Tsallis bridge pb}
\inf_{\pi \in \Pi(\mu_1, \mu_2)} \ep D_q(\pi, \mathcal{K}_q), 
\end{equation}
where $\mathcal{K}_q$ is defined by 
\begin{equation*}
\mathcal{K}_q = \exp_{q} \left(- \dfrac{c(x,y)}{\ep} \right)  \mu_1 \otimes \mu_2 = \left[ 1 + (q-1) \left( \dfrac{c(x,y)}{\ep} \right) \right]^{-\frac{1}{q-1}} \mu_1 \otimes \mu_2. 
\end{equation*}
Using the non-additivity property of the $q$-logarithmic function
\begin{equation*}
\log_{2-q} (xy) = \log_{2-q} x + \log_{2-q} y + (q-1) \log_{2-q} x \log_{2-q} y
\end{equation*}
and $\log_{2-q} \left( \exp_{q}(-x)^{-1} \right) = x$, the minimization problem \eqref{Tsallis bridge pb} corresponds to the Tsallis entropic regularized optimal transport problem \eqref{Tsallis reg OT} as follows: 
\begin{equation*}
\mathrm{OT}_{q,\ep} = \inf_{\pi \in \Pi(\mu_1, \mu_2)} \left\{ \ep D_q(\pi, \mathcal{K}_q) - (q-1) \int_{X_1 \times X_2} c \log_{2-q} \left( \dfrac{d \pi}{d (\mu_1 \otimes \mu_2 )}\right) d\pi \right\}. 
\end{equation*}
We note that the equivalence between the Tsallis entropic regularized optimal transport problem \eqref{Tsallis reg OT} and the minimization problem \eqref{Tsallis bridge pb} involves an interaction term, whereas the KL entropic regularized optimal transport \eqref{KL reg OT} is equivalent to the Schr\"{o}dinger bridge problem \eqref{Schrodinger bridge pb}. Establishing a quantitative estimate between \eqref{Tsallis reg OT} and \eqref{Tsallis bridge pb} would be an interesting direction for future research.

\section*{Data availability statement}
Data sharing not applicable to this article as no datasets were generated or analysed during the current study.

\section*{Acknowledgements}
The authors would like to thank Dr.~Jun Okamoto (Kyoto University) for the helpful discussions. The first author is partially supported by JSPS Research Activity Start-up (22K20336).
The second author is partially supported by JSPS Grant-in-Aid for Early-Career Scientists (21K13822), JSPS Grant-in-Aid for Transformative Research Areas (A) (22A201), and WPI-ASHBi at Kyoto University. 


\bibliography{Tsallis}
\bibliographystyle{siam}

\bigskip

\noindent
\textsc{
Faculty of Advanced Science and Technology, Kumamoto University, Kumamoto 860-8555, Japan}\\
\noindent
{\em Electronic mail address:}
suguro@kumamoto-u.ac.jp

\bigskip

\noindent
\textsc{ 
Mathematical Science Center for Co-creative Society, Tohoku University, 
Sendai 980-0845, Japan} \\
\noindent
{\em Electronic mail address:}
toshiaki.yachimura.a4@tohoku.ac.jp

\end{document}